\documentclass[a4paper,11pt]{article}

\usepackage{graphicx}
\usepackage{bm}
\usepackage[percent]{overpic}

\usepackage{authblk}
\usepackage{amscd}
\usepackage{amsfonts}
\usepackage[numbers]{natbib}
\usepackage{amsthm}
\usepackage{amssymb}
\usepackage{amsmath}
\usepackage{algpseudocode}
\usepackage{algorithm}
\usepackage{cancel}
\usepackage{chngcntr}
\usepackage{caption}
\usepackage{enumerate}
\usepackage{epsfig}
\usepackage{epic}
\usepackage{eepic}
\usepackage{enumitem}
\usepackage{etoolbox}
\usepackage{fancyhdr}
\usepackage{fullpage}
\usepackage[a4paper, total={6in, 8in}]{geometry}
\usepackage{graphicx}
\usepackage{here}
\usepackage[utf8]{inputenc}
\usepackage{lscape}
\usepackage{listings}
\usepackage{longtable}
\usepackage{natbib}
\usepackage{nccmath}
\usepackage{parskip}
\usepackage{rotating}
\usepackage{subfigure}
\usepackage{tabularx}
\usepackage{textcomp}
\usepackage{tikz}
\usepackage{titlesec}
\usepackage{threeparttable}
\usepackage{varwidth}
\usepackage{framed,multirow}
\graphicspath{ {./Figures_FEM/} }
\newtheorem{pro}{Problem}
\usepackage[percent]{overpic}
\usepackage{url}

\usepackage{xcolor}

\title{Finite element discretization of a biological network formation system: a preliminary study}

\author[1]{Clarissa Astuto}
\author[1,2]{Daniele Boffi}
\author[1,2]{Fabio Credali}

\affil[1]{King Abdullah University of Science and Technology (KAUST), 4700, Thuwal, Saudi Arabia}
\affil[2]{Department of Mathematics ''F. Casorati'', University of Pavia, Pavia, Italy}

\begin{document}
\maketitle

\begin{abstract}
A finite element discretization is developed for the Cai-Hu model, describing the formation of biological networks. The model consists of a non linear elliptic equation for the pressure $p$ and a non linear reaction-diffusion equation for the conductivity tensor $\mathbb{C}$. The problem requires high resolution due to the presence of multiple scales, the stiffness in all its components and the non linearities. We propose a low order finite element discretization in space coupled with a semi-implicit time advancing scheme. The code is {verified} with several numerical tests performed with various choices for the parameters involved in the system. In absence of the exact solution, we apply Richardson extrapolation technique to estimate the order of the method.
\end{abstract}

\section{Introduction}
The branching structure of biological transport networks, such as leaf venations and blood flow, has been the object of extensive research due to their widespread existence in living organisms and natural phenomena,{ such as, leaf venation in plants \cite{malinowski2013understanding},
	vascular pattern formation \cite{sedmera2011function}, mammalian circulatory systems and
	neural networks transporting electric charge \cite{eichmann2005guidance, michel1995morphogenesis}}.

It has been suggested that the optimization of energy consumption, as a result of natural selection, governs the branching structure. In 2013, Hu and Cai \cite{hu2013adaptation} proposed a global energy functional approach that {takes into consideration a metabolic cost to sustain the network. This consumption is also influenced by the the local gradient of pressure that can modify the diameter of the vessels as an adaptive response to the stress \cite{hu2012blood}.}


{The considered equations have been recently analyzed, from the numerical point of view, making use of finite difference schemes \cite{astuto2023asymmetry,astuto2022comparison}. Finite element methods (FEM) (see, e.g., \cite{ciarlet2002finite}) are appealing for their flexibility and {are} widely used by the mathematical community. In this work we present a low order finite element method for solving the biological network PDE system in tensor form and capture the steady state solution.} The approximation of biological network formation has been already studied {by means} of mixed finite elements \cite{bellomo2017active,marko_perthame_schlo} in the case of the corresponding vector model. A comparison between the two different versions, tensor and vector models, has been showed in \cite{astuto2022comparison}. 

Our code is based on the finite elements library FreeFEM \cite{PironneauFREEFEM,MR3043640}, which proposes a large variety of triangular
finite elements to solve partial
differential equations. {FreeFEM has been already used for numerical computations of stationary states of fast-rotating Bose–Einstein condensates \cite{DANAILA20106946}, and for the minimization of energy functionals associated with Ginzburg-Landau and Schr\"odinger models \cite{raza2010approximating,raza2009energy}}.

The paper is organized as follows: in Section~\ref{section_model}, we introduce the tensor PDE model and we show the derivation of the weak formulation. In Section~\ref{section_discretization}, we describe the discrete version of the problem by coupling finite elements in space with a semi-implicit time advancing scheme based on Backward Euler.
In Section~\ref{section_results}, we report and discuss accuracy tests and results of the numerical simulations. Finally, in Section~\ref{section_conclusion}, we summarize the results and draw some conclusions.

\section{Mathematical model}
\label{section_model}
We consider an elliptic-parabolic system of partial differential equations for the conductivity tensor, $\mathbb{C}\in \mathbb{R}^{2\times 2}$, 
and a pressure scalar force, $p\in \mathbb{R}$. The model, introduced in \cite{hu2013adaptation}, reads
\begin{align}
	\label{eq_darcy_p_fin}
	-\nabla\cdot \left( (r\mathbb{I}+\mathbb{C}) \nabla p \right) &= S \qquad {\rm in }\, \Omega \\
	\label{eq_reaction_diff_fin}
	\frac{\partial\mathbb{C}}{\partial t} - D^2\Delta \mathbb{C} -c^2\nabla p \otimes \nabla p + \alpha \|\mathbb{C}\|^{\gamma-2} \mathbb{C} & = 0  \qquad {\rm in }\, \Omega.
\end{align}
The term $S = S(\vec{x})$ denotes the distribution of sources and sinks, which is a known datum. The function $r :\Omega \to \mathbb{R}^+$, with $r(x) \geq r_0 > 0$, describes the isotropic background permeability of the medium.
In Eq.~\eqref{eq_reaction_diff_fin} the diffusion coefficient ${D>0}$ controls the random effects in the transportation medium, while the activation parameter $c > 0$ describes the tendency of the network to align with the pressure gradient.
The reaction term $\alpha \|\mathbb{C}\|^{\gamma-2} \mathbb{C}$, where $\|\cdot\|$ denotes the Frobenius norm, models the metabolic cost of maintaining the network structure, with metabolic coefficient $\alpha>0$ and metabolic exponent $\gamma>0$. For instance, for leaf venation in plants, $1/2\leq \gamma \leq 1$  \cite{hu2013adaptation, hu2013optimization}. {For this range of values, there is no proof of convergence to a unique steady state, while it has been proved for $\gamma>1$ in \cite{marko_perthame_2}. Both cases have been numerically investigated in \cite{astuto2022comparison}.}


We consider Eqs.(\ref{eq_darcy_p_fin}--\ref{eq_reaction_diff_fin})
in a bounded domain $\Omega \subset \mathbb{R}^2$, with smooth boundary $\partial \Omega$.
We choose homogeneous Neumann boundary conditions for $p$ and $\mathbb{C}$  
\begin{equation}
	\label{eq_bc}
	\nabla p(t,\vec{x}) \cdot \nu = 0, \quad \nabla \mathbb{C}(t,\vec{x})\cdot \nu = 0 , \quad \text{for } \vec{x}\in \partial \Omega, \, t\geq 0,
\end{equation}
where $\nu$ denotes the outer normal vector to $\partial \Omega$. With the homogeneous Neumann boundary condition for
the pressure $p$, we impose the global mass balance
\begin{equation}
	\label{eq_Source}
	\int_\Omega S(\vec{x})\,d\vec{x} = 0
\end{equation}
to ensure solvability of Eq.~\eqref{eq_darcy_p_fin}.
{We choose Neumann boundary conditions for $\mathbb{C}$ because in \cite{astuto2022comparison} we observed that homogeneous Dirichlet boundary conditions lead to the formation of boundary layers in the solutions. On the contrary, homogeneous Neumann boundary conditions seem able to bypass this problem. Let us also notice that the same type of conditions was selected for the numerical simulations carried out in \cite{hu2019optimization}.
	Homogeneous Neumann boundary conditions for the pressure $p$ are the most natural choice since we impose zero flux of the transported material through $\partial\Omega$, where the material flux is given by ${\vec q = (r\mathbb{I}+\mathbb{C} )\nabla p}$.}

At the end, the system is closed by defining an initial condition $\mathbb{C}_0 \geq 0$ for the conductivity $\mathbb{C}$,
\begin{equation}
	\label{eq_ic}
	\mathbb{C}(t=0,\vec{x}) = \mathbb{C}_0(\vec{x}) \quad \text{ in } \Omega.
\end{equation}

The system is derived as gradient flow of an energy functional, consisting of a diffusive term, an activation term and a metabolic cost term. The expression for the functional $\mathcal{E}[\mathbb{C}]$ reads
\begin{equation} \label{energy_T}
	\mathcal{E}[\mathbb{C}] :=
	\int_{\Omega } \left( \frac{D^2}{2} |\nabla \mathbb{C}|^2 + c^2 \nabla p[\mathbb{C}] \cdot \mathbb{P}[\mathbb{C}] \nabla p[\mathbb{C}]  + \frac{\alpha}{\gamma}|\mathbb{C}|^\gamma \right) \mathrm{d} \vec{x},
\end{equation}
where $\mathbb{P}[\mathbb{C}] = r\mathbb{I}+\mathbb{C}$, and, in Section~\ref{section_results}, {we show that our numerical method provides energy decay, in agreement with the mathematical model.}  

Before presenting the variational formulation of Eqs.(\ref{eq_darcy_p_fin}--\ref{eq_reaction_diff_fin}), we discuss the functional analysis setting in which we are going to work. 

The variable $p$ is the solution of the diffusion equation with Neumann boundary conditions, consequently we impose a null-mean condition that guarantees uniqueness of the solution. Therefore, we introduce the space $Q$ defined as
\begin{equation}
	Q = \biggl\{ q\in H^1(\Omega) : \int_\Omega q\,{ \rm d}\Omega = 0\biggr\}.
\end{equation}
For the conductivity tensor $\mathbb{C}$, we introduce the following tensor-valued space, which takes into account the symmetry
\begin{equation}
	\mathbf{V} = \biggl\{ \mathbb{B}\in[H^1(\Omega)]^{2\times2}:\mathbb{B}=\mathbb{B}^\top \biggr\}.
\end{equation}
Multiplying Eqs.(\ref{eq_darcy_p_fin}--\ref{eq_reaction_diff_fin}) by test functions $q\in Q$ and $\mathbb{B}\in\mathbf{V}$, respectively, and integrating, we obtain

\begin{align}
	&-\int_\Omega \nabla\cdot \left( (r\mathbb{I}+\mathbb{C}) \nabla p \right)\, q\,{ \rm d}\Omega  = \int_\Omega S\, q\,{ \rm d}\Omega && \forall q\in Q \\
	&\int_\Omega \frac{\partial \mathbb{C}}{\partial t} \, \mathbb{B}\,{ \rm d}\Omega - D^2\int_\Omega (\Delta \mathbb{C}) \, \mathbb{B}\,{ \rm d}\Omega\\
	& \notag \hspace{1.5cm} - c^2\int_\Omega (\nabla p \otimes \nabla p) \, \mathbb{B}\,{ \rm d}\Omega + \alpha \int_\Omega \| \mathbb{C} \|^{\gamma - 2} \mathbb{C} \, \mathbb{B}\,{ \rm d}\Omega  = 0 && \forall \mathbb{B}\in \mathbf{V} .
\end{align}
Integrating by parts and taking into account the boundary conditions, we have
\begin{align}
	&\int_\Omega (r\mathbb{I}+\mathbb{C}) \nabla p \cdot \nabla q\,{ \rm d}\Omega  = \int_\Omega S\, q\,{ \rm d}\Omega && \forall q\in Q \\
	&\int_\Omega \frac{\partial \mathbb{C}}{\partial t} \, \mathbb{B}\,{ \rm d}\Omega + D^2\int_\Omega \nabla \mathbb{C} : \nabla\mathbb{B}\,{ \rm d}\Omega\\
	&\notag\hspace{1.5cm}- c^2\int_\Omega (\nabla p \otimes \nabla p) \, \mathbb{B}\,{ \rm d}\Omega + \alpha \int_\Omega \| \mathbb{C} \|^{\gamma - 2} \mathbb{C} \, \mathbb{B}\,{ \rm d}\Omega  = 0 && \forall \mathbb{B}\in \mathbf{V} .
\end{align}
Finally, denoting by $(\cdot,\cdot)$ the scalar product in $L^2(\Omega)$, our problem in variational formulation reads as follows.

\begin{pro}\label{pro:variational}
	Given $S \in L^2(\Omega)$ and $\mathbb{C}_0\in\mathbf{V}$, such that $\nabla \mathbb{C}^0(t,\vec{x})\cdot \nu = 0$, find $p(t)\in Q$ and $\mathbb{C}(t) \in \mathbf{V}$ such that, for almost every $t\in(0,T)$, it holds
	\begin{align}
		&\big( (r\mathbb{I}+\mathbb{C}) \nabla p , \nabla q \big)  = \big( S,q \big) && \forall q\in Q \\
		&\bigg(\frac{\partial \mathbb{C}(t)}{\partial t},\mathbb{B}\bigg) + D^2\big(\nabla \mathbb{C}(t),\nabla\mathbb{B}\big)\\
		&\notag\hspace{1.5cm}- c^2\big(\nabla p(t) \otimes \nabla p(t),\mathbb{B}\big) + \alpha \big( \| \mathbb{C}(t) \|^{\gamma - 2} \mathbb{C}(t) , \mathbb{B}\big) = 0 && \forall \mathbb{B}\in \mathbf{V}\\
		& \nabla p(t,\vec{x}) \cdot \nu = 0 && \text{on } \partial\Omega\\
		& \nabla \mathbb{C}(t,\vec{x})\cdot \nu = 0 && \text{on } \partial\Omega\\
		& \mathbb{C}(0,\vec{x}) = \mathbb{C}_0(\vec{x}) && \text{in } \Omega.
	\end{align}
\end{pro}

\section{Discretization}
\label{section_discretization}
We discretize Problem~\ref{pro:variational} in space by finite elements. Let us consider a mesh $\mathcal{T}_h$ of $\Omega$ with mesh size $h$, so that we can introduce two finite dimensional subspaces $Q_h\subset Q$ and $\mathbf{V}_h\subset\mathbf{V}$. The discrete problems is the following one.

\begin{pro}\label{pro:variational_2}
	Given $S \in L^2(\Omega)$ and $\mathbb{C}_{0,h}\in\mathbf{V}_h$, such that $\nabla \mathbb{C}_h^0(t,\vec{x})\cdot \nu = 0$, find $p_h(t)\in Q_h$ and $\mathbb{C}_h(t) \in \mathbf{V}_h$ such that, for almost every $t\in(0,T)$, it holds
	\begin{align}
		&\label{eq:p_disc}\big( (r\mathbb{I}+\mathbb{C}_h) \nabla p_h,\nabla q_h\big)  = \big( S, q_h\big) && \forall q_h\in Q_h \\
		&\label{eq:C_disc}\bigg( \frac{\partial \mathbb{C}_h(t)}{\partial t} ,\mathbb{B}_h\bigg) + D^2\big(\nabla \mathbb{C}_h(t) , \nabla\mathbb{B}_h\big) - c^2\big(\nabla p_h(t) \otimes \nabla p_h(t),\mathbb{B}_h\big)\\
		&\notag\hspace{4.5cm} + \alpha \big( \| \mathbb{C}_h(t) \|^{\gamma - 2} \mathbb{C}_h(t) , \mathbb{B}_h\big)  = 0 && \forall \mathbb{B}_h\in \mathbf{V}_h\\
		& \nabla p_h(t,\vec{x}) \cdot \nu = 0 && \text{on } \partial\Omega\\
		& \nabla \mathbb{C}_h(t,\vec{x})\cdot \nu = 0 && \text{on } \partial\Omega\\
		& \mathbb{C}_h(0,\vec{x}) = \mathbb{C}_{0,h}(\vec{x}) && \text{in } \Omega.
	\end{align}
\end{pro}

We now discuss the discretization in time. We consider a uniform partition of the time interval $[0,T]$ in $N\in\mathbb{N}$ parts of length $\Delta t$ denoting the nodes by $t^n = n\Delta t$ and $\mathbb{C}_h^n \approx \mathbb{C}_h(t^n)$. A first order semi-implicit approximation of the time derivative of $\mathbb{C}_h$ is obtained with Implicit Euler scheme
\begin{equation*}
	\frac{\partial \mathbb{C}_h(t_{n+1})}{\partial t}\approx\frac{\mathbb{C}_h^{n+1} -\mathbb{C}_h^{n}}{\Delta t}
\end{equation*}
so that Eq.~\eqref{eq:C_disc} can be rewritten as 
\begin{align*}
	\left(\, \mathbb{C}_h^{n+1}, \mathbb{B}_h \,\right) + \Delta t\, D^2 \left(\, \nabla \mathbb{C}_h^{n+1}, \nabla \mathbb{B}_h \,\right) &- \Delta t\, c^2 \left(\, \nabla p_h^{n} \otimes \nabla p_h^{n}, \mathbb{B}_h\,\right)\\&+ \Delta t\, \alpha \left(\, \left(\| \mathbb{C}_h^{n} \| + \varepsilon \right)^{\gamma - 2} \mathbb{C}_h^{n+1},\mathbb{B}_h \right)
	= \left(\, \mathbb{C}_h^{n}, \mathbb{B}_h \,\right),
\end{align*}
where we add a \textit{regularization parameter} $\varepsilon > 0$, to prevent the instability coming from the division by zero {in the metabolic term} when $\gamma < 1$ {and $\mathbb{C}$ is close to zero in the region without ramifications} (see \cite{astuto2022comparison,astuto2023asymmetry} for more details). In order to solve Problem~\ref{pro:variational_2}, we adopt a partitioned scheme; in particular, the nonlinear term $\left(\, \nabla p_h^{n} \otimes \nabla p_h^{n}, \mathbb{B}_h\,\right)$ is computed by means of $p_h^n$, which is obtained by solving 
\[
\left(\, \left(r\mathbb{I}+\mathbb{C}_h^n\right) \nabla p_h^n ,\nabla q_h \,\right)  = \left(\, S , q_h \right)
\]
at the previous time step. Moreover, we remark that the nonlinear metabolic term $$\left(\, \left(\| \mathbb{C}_h^{n} \| + \varepsilon\right)^{\gamma - 2} \mathbb{C}_h^{n+1},\mathbb{B}_h \right)$$ is treated explicitly.

In this setting, our time advancing scheme is of order one. Anyway, a second order scheme cannot be easily obtained because of the nonlinearity of the system. Fully second order accuracy could be recovered by adopting, for example, a predicted pressure at time $t^{n+1/2}$, as suggested in \cite{CiCP-31-707}, or by treating the metabolic term by means of a semi-implicit IMEX approach (see, e.g., \cite{IMEX}). 

\section{Results}
\label{section_results}
In this section we perform several simulations with different choices of the parameters involved in the discrete model. In particular, we see that the results agree with the previous work obtained with finite difference schemes \cite{astuto2023asymmetry,astuto2022comparison}. 

The domain we consider is $\Omega = [0,1]\times [0,1]$ and the initial conditions and source term are defined in Eqs.(\ref{eq_ic_C}--\ref{eq_S_term}). In Table~\ref{tab_qualit} we collect the values of all the parameters related to our tests. 
\begin{table}
	\centering
	\begin{tabular}{|c|c|c|c|c|c|c|c|}  \hline
		& $\alpha$ & $c^2$ & $D^2$ & $\varepsilon$ & $\gamma$ & $r$ & $T$ \\ \hline 
		\textsc{Test-A1}: & 0.75 & 1 & $10^{-1}$ & $10^{-1}$ & 1.25 & $10^{-1}$ & 1 \\ \hline
		\textsc{Test-A2}: &  0.75 & 1 & $10^{-1}$ & $10^{-1}$ & 0.75 & $10^{-1}$ & 1 \\ \hline 
		\textsc{Test-A3}: &  0.75 & 1 & $10^{-4}$ & $10^{-1}$ & 0.75 & $10^{-1}$ & 1 \\ \hline \hline
		\textsc{Test-Al1}: &  1.5 & 25 & $10^{-3}$ & $10^{-3}$ & 0.75 & $10^{-3}$ & 10 \\ \hline
		\textsc{Test-Al2}: &  1.5 & 25 & $5\cdot 10^{-4}$ & $10^{-3}$ & 0.75 & $10^{-3}$ & 10 \\ \hline
		\hline              
		\textsc{Test-DD1}: &  0.75 & 25 & $2\cdot 10^{-3}$ & $10^{-3}$ & 0.75 & $10^{-3}$ & 10 \\ \hline
		\textsc{Test-DD2}: &  0.75 & 25 & $ 10^{-3}$ & $10^{-3}$ & 0.75 & $10^{-3}$ & 10 \\ \hline
		\textsc{Test-DD3}: &  0.75 & 25 & $5\cdot 10^{-4}$ & $10^{-3}$ & 0.75 & $10^{-3}$ & 10 \\ \hline \hline         \textsc{Test-}r1: &  0.75 & 25 & $10^{-3}$ & $10^{-3}$ & 0.75 & $10^{-1}$ & 10 \\ \hline
		\textsc{Test-}r2: &  0.75 & 25 & $10^{-3}$ & $10^{-3}$ & 0.75 & $10^{-2}$ & 10 \\ \hline 
	\end{tabular}
	\caption{\textit{Parameters involved in the tests shown in Section~\ref{section_results}. The first three rows are related to the accuracy tests: their results are summarized in Table~\ref{table_accuracy}. 
	}}
	\label{tab_qualit}
\end{table}

We define the initial condition for each component of the tensor $\mathbb{C}_0(\vec{x})$ and the source function $S(\vec{x})$ as follows:
\begin{eqnarray} \label{eq_ic_C}
	& \mathbb{C}_0^{(1,1)}(\vec{x
	}) = 1, \quad \mathbb{C}_0^{(1,2)}(\vec{x
	}) = 0, \quad \mathbb{C}_0^{(2,2)}(\vec{x
	}) = 1,  \\ \label{eq_S_term}
	& S(\vec{x}) = E(\vec{x}) - \overline{E}, \quad
	E(\vec{x}) = \exp(-\sigma(\vec{x}-\vec{x}_0)^2),\\
	& \sigma = 500, \quad \vec{x}_0 = (0.25,0.25),
\end{eqnarray}
where $\overline{E}$ is the average of $E$ over the domain $\Omega$. {We remark that the initial conditions and sink-source function $S$ are identical to those already used in \cite{astuto2023asymmetry}. Moreover, \textsc{Test-DD3} is reproduced from the same paper.}

The discretization of the problem is performed by considering continuous piecewise linear finite elements for both $p$ and $\mathbb{C}$ on uniform triangulations $\mathcal{T}_h$ of $\Omega$; in particular, we have
\begin{equation}
	\begin{aligned}
		& Q_h = \{ q\in H^1(\Omega):\,q_{|K}\in\mathcal{P}_1(K)\,\,\forall K\in\mathcal{T}_h\}\\
		& \mathbf{V}_h = \{ \mathbb{B}\in [H^1(\Omega)]^{2\times 2}:\,\mathbb{B}_{|K}\in[\mathcal{P}_1(K)]^{2\times2}\,\,\forall K\in\mathcal{T}_h \}.
	\end{aligned}
\end{equation}

\subsection{Accuracy tests}
\label{section_qualit_test}
First, we check the accuracy in time of the adopted scheme. In Table~\ref{table_accuracy}, we see the error for the norm of the conductivity tensor, calculated with Richardson extrapolation (see, e.g.,~\cite{richardson1911ix}) by refining the mesh size $h$. We remark that the time step $dt$ depends on $h$, since we have $dt=h$. The chosen parameters are collected in \textsc{Test-A1}--\textsc{Test-A3}.

\begin{table}
	\centering
	\begin{tabular}{||c||c||c||} \hline \hline
		$h$ & error & order \\ \hline \hline
		1/20 & -- & -- \\ \hline \hline
		1/40 & 0.0117878 & -- \\ \hline \hline
		1/80 & 0.00612581 & 0.94 \\ \hline \hline
		1/160 & 0.00312414 & 0.97 \\ \hline \hline
		1/320 & 0.00157781 & 0.99\\ \hline \hline
	\end{tabular}
	\begin{tabular}{||c||c||c||} \hline \hline
		$h$ & error & order \\ \hline \hline
		1/20 & -- & -- \\ \hline \hline
		1/40 & 0.0136681  & -- \\ \hline \hline
		1/80 & 0.00721071 &  0.92 \\ \hline \hline
		1/160 & 0.00370715 & 0.96 \\ \hline \hline
		1/320 & 0.00188007 & 0.98 \\ \hline \hline
	\end{tabular}
	\begin{tabular}{||c||c||c||} \hline \hline
		$h$ & error & order \\ \hline \hline
		1/20 & -- & -- \\ \hline \hline
		1/40 & 0.0136681  & -- \\ \hline \hline
		1/80 & 0.00721071 &  0.92 \\ \hline \hline
		1/160 & 0.00370715 & 0.96 \\ \hline \hline
		1/320 & 0.00188007 & 0.98 \\ \hline \hline
	\end{tabular}
	\caption{\textit{{Accuracy study by computing the relative error in $L^2-$norm between two solutions with two different mesh sizes. In order to compute the error for $\|\mathbb{C}\|$ on $\mathcal{T}_h$, we consider as reference solution $\|\mathbb{C}\|$ computed in the previous test on $\mathcal{T}_{h/2}$. The problem has been solved with parameters defined in  \textsc{Test-A1}--\textsc{Test-A3}.}}}
	\label{table_accuracy}
\end{table}

\subsection{Numerical tests}
In this subsection, we show and analyze the results of our tests in terms of $\|\mathbb{C}\|$: in particular, we analyze its snapshots when the steady state is reached. {When this happens, the energy functional remains constant in time. {Another} possible indicator may be the numerical time derivative of the tensor, that should be approximately close to zero.} 

Again, the parameters of the tests are collected in Table~\ref{tab_qualit}, the initial conditions and source term in Eqs.(\ref{eq_ic_C}--\ref{eq_S_term}). Each test has been performed by considering a triangulation with size $h=1/600$ and $dt=0.01$.
\\
\begin{figure}
	\vspace{5mm}
	\centering
	\begin{minipage}{.3\textwidth}
		\centering
		\begin{overpic}[abs,width=0.9\textwidth,unit=1mm,scale=.25]{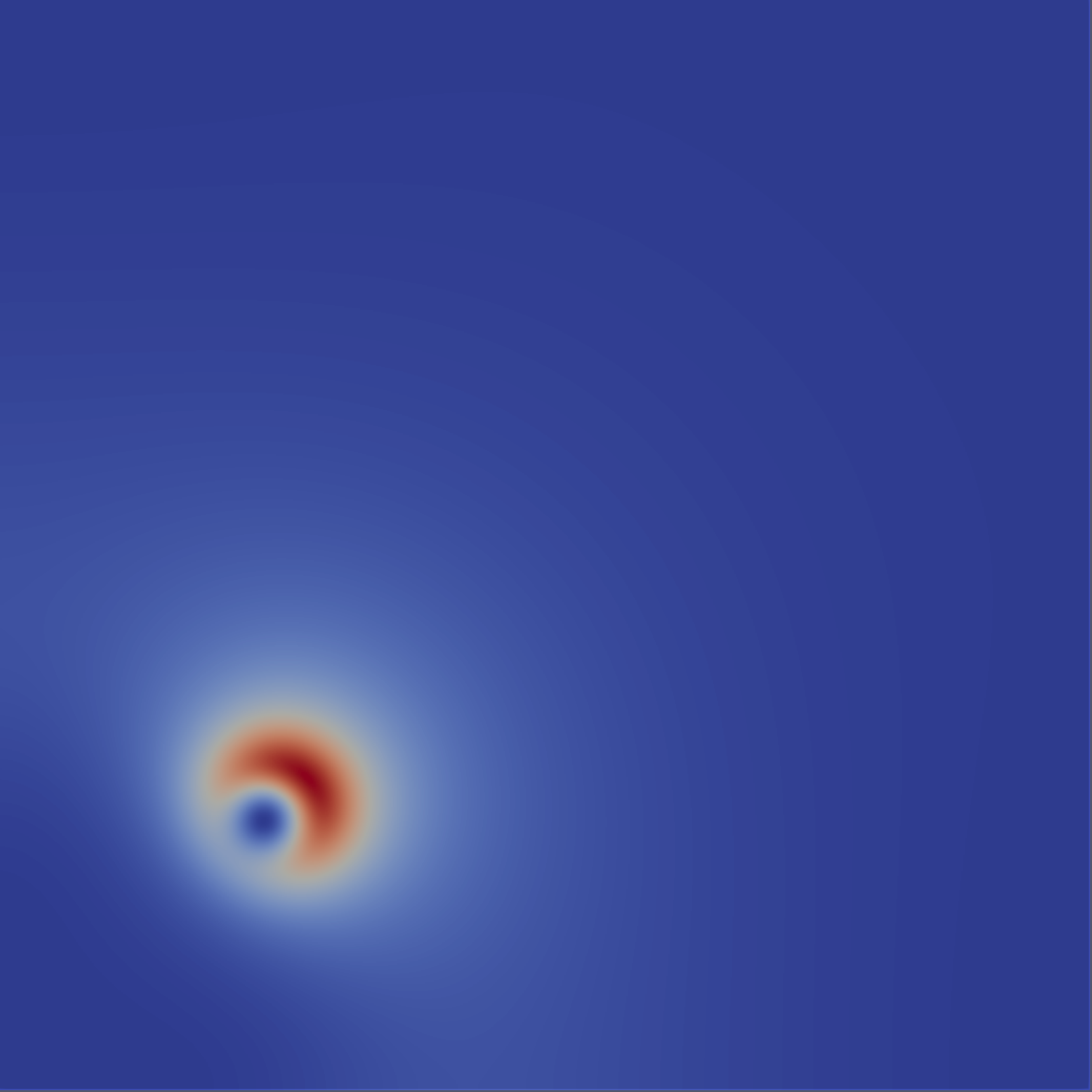}
			\put(0.5,0.1){\textcolor{white}{$\rightarrow\,x$}}
			\put(0.01,2){\textcolor{white}{$\uparrow$}}
			\put(0.01,7){\textcolor{white}{$y$}}
			\put(9,43){\textsc{Test-}r1: $\|\mathbb{C}\|$}
		\end{overpic}
	\end{minipage}
	\begin{minipage}{.3\textwidth}   
		\centering
		\begin{overpic}[abs,width=0.9\textwidth,unit=1mm,scale=.25]{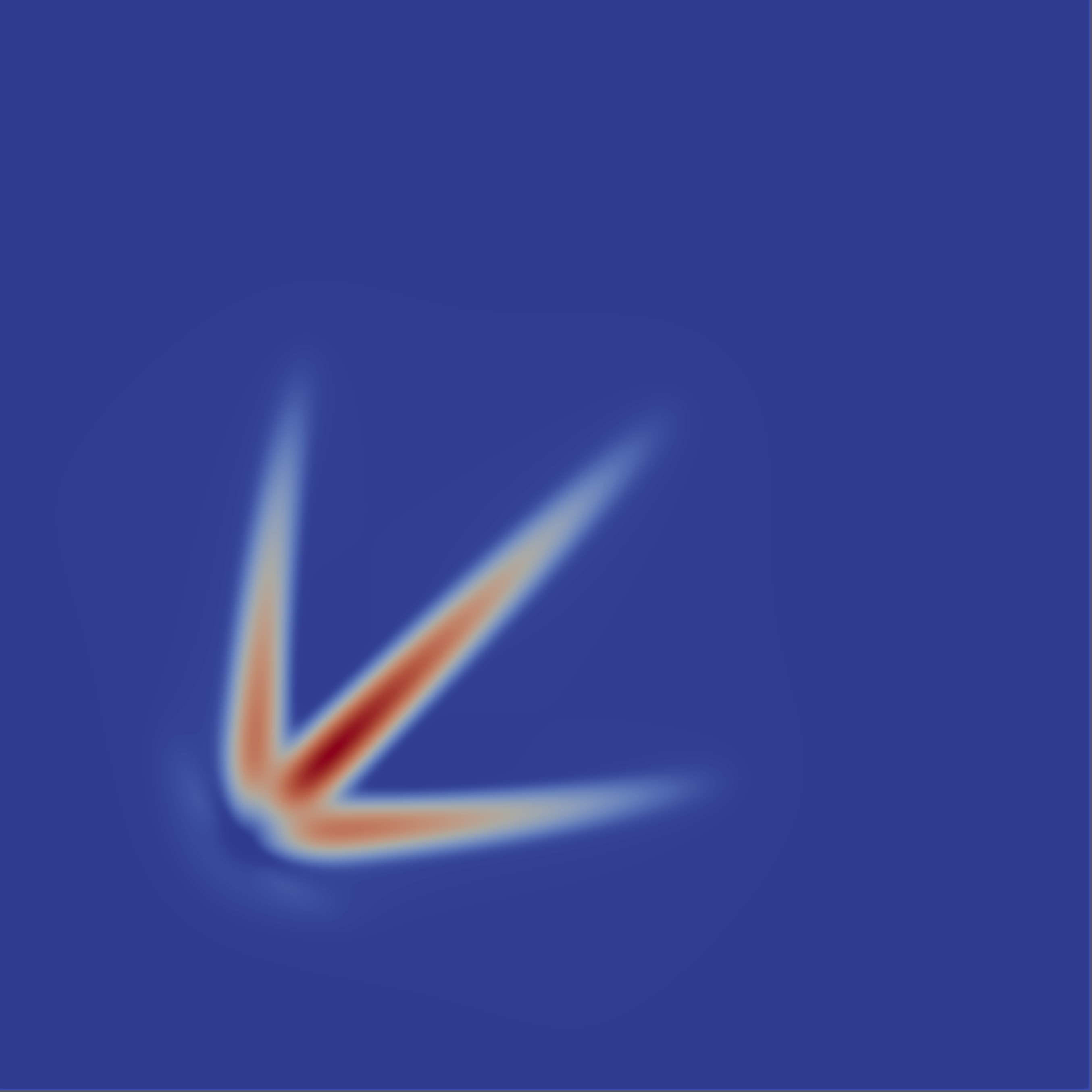}
			\put(0.5,0.1){\textcolor{white}{$\rightarrow\,x$}}
			\put(0.01,2){\textcolor{white}{$\uparrow$}}
			\put(0.01,7){\textcolor{white}{$y$}}
			\put(9,43){\textsc{Test-}r2: $\|\mathbb{C}\|$}
		\end{overpic}
	\end{minipage}
	\begin{minipage}{.3\textwidth} 
		\centering
		\begin{overpic}[abs,width=0.9\textwidth,unit=1mm,scale=.25]{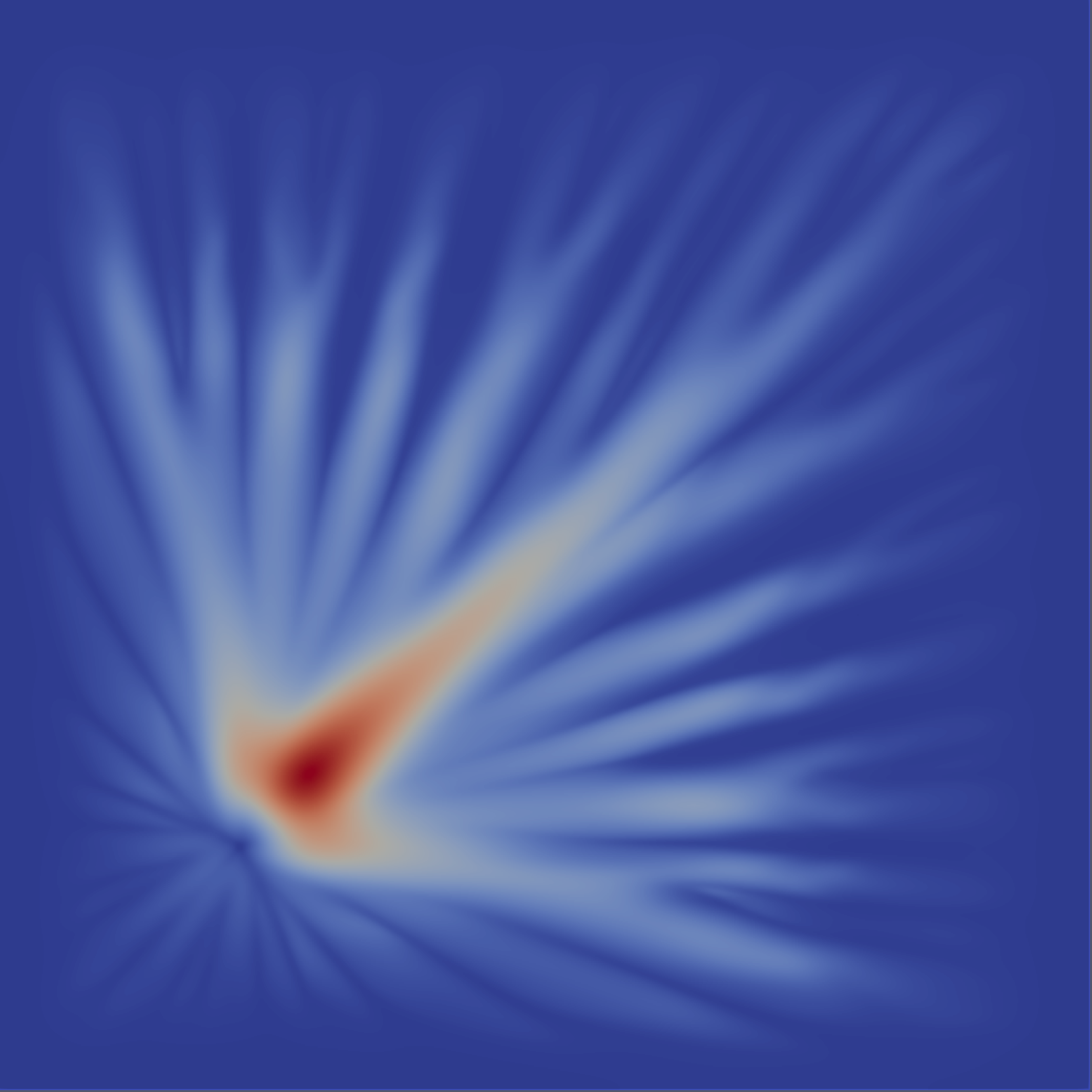}
			\put(0.5,0.1){\textcolor{white}{$\rightarrow\,x$}}
			\put(0.01,2){\textcolor{white}{$\uparrow$}}
			\put(0.01,7){\textcolor{white}{$y$}}
			\put(7,43){\textsc{Test-DD1}: $\|\mathbb{C}\|$}
		\end{overpic}
	\end{minipage}
	\begin{minipage}{.3\textwidth} \centering\includegraphics[width=1.0\textwidth]{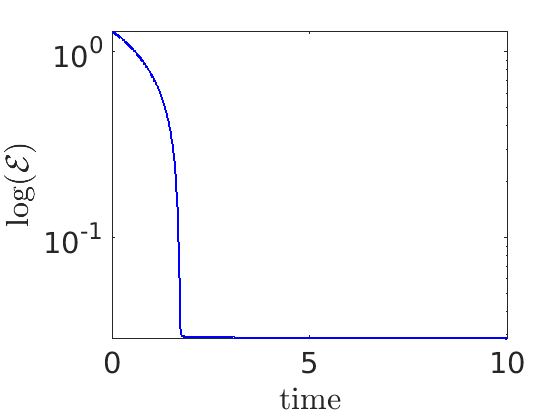}
	\end{minipage}
	\begin{minipage}{.3\textwidth} \centering   \includegraphics[width=1.0\textwidth]{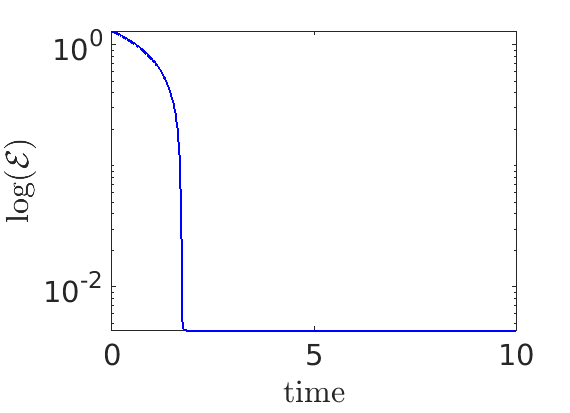}
	\end{minipage}
	\begin{minipage}{.3\textwidth}  \centering  \includegraphics[width=1.0\textwidth]{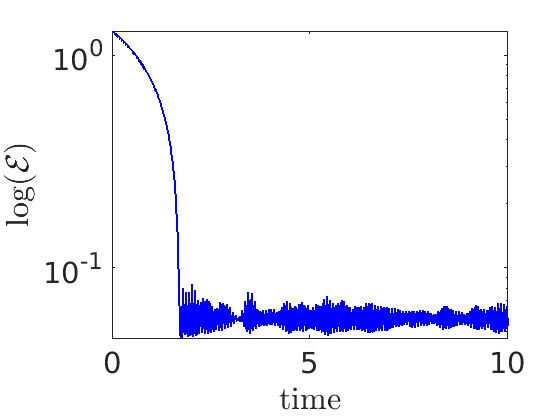}
	\end{minipage}
	\caption{\textit{In this figure we show three results considering different values of the background permeability $r$. From left to right we have the results of} \textsc{Test-}r1 \textit{with $r = 10^{-1}$,} \textsc{Test-}r2 \textit{with $r = 10^{-2}$ and} \textsc{Test-DD1} \textit{with $r = 10^{-3}$. We plot two different quantities for each test: in the first row, the norm of the solution $\mathbb{C}$ at final time; in the second row, the time dependent plot of the energy $\mathcal{E}[\mathbb{C}]$.}}
	\label{fig_testr}
\end{figure}
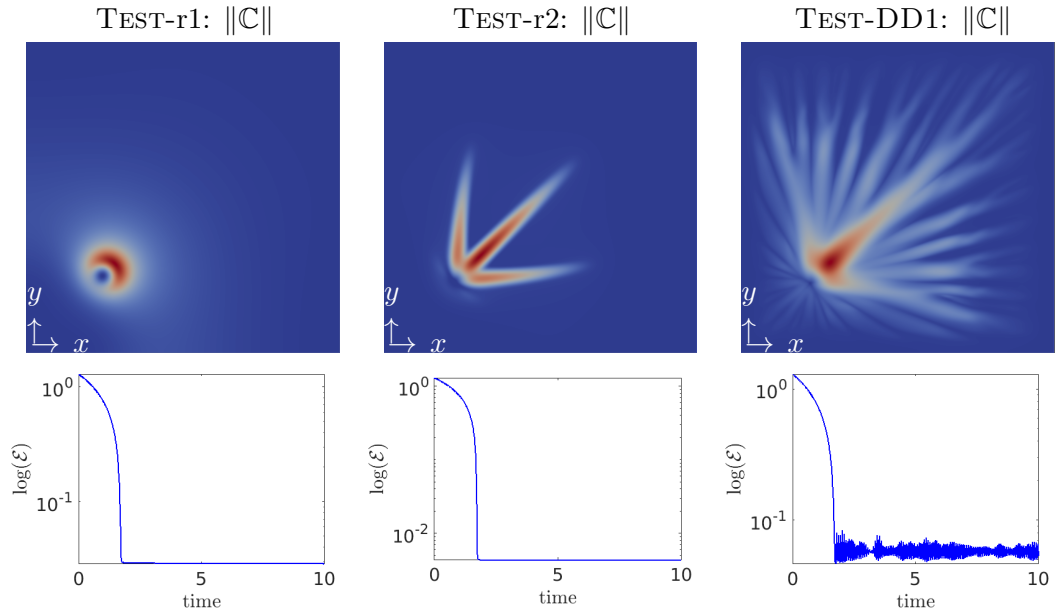
In Fig.~\ref{fig_testr} we show the results with different values of the background permeability $r$. The tests we consider are \textsc{Test-}r1, \textsc{Test-}r2 and \textsc{Test-DD1}. In terms of $\|\mathbb{C}\|$, the most interesting case is \textsc{Test-DD1}, with $r = 10^{-3}$, because {a more detailed pattern appears}, in agreement with the results in \cite{astuto2023asymmetry,astuto2022comparison}. When $r \to 0$, the energy plots show some oscillations while approaching the steady state. Moreover, as showed in \cite{astuto2023asymmetry}, the iteration matrix for the pressure becomes ill conditioned when we decrease the parameter $r$, and this was causing the deterioration of the symmetry of the solution. In this case, not only we see the lack of symmetry, but also the non monotonic behaviour of the energy in time.

In Fig.~\ref{fig_testDD}, we show the results with different values of the diffusivity $D^2$. The tests we consider are \textsc{Test-DD1} with $D^2 = 2\cdot 10^{-3}$, \textsc{Test-DD2} with $D^2 = 10^{-3}$ and \textsc{Test-DD3} $D^2 = 5\cdot 10^{-4}$. From the plot of $\|\mathbb{C}\|$, we can see that the ramifications become thinner and thinner when decreasing the diffusivity, while the plots of the energy show some oscillations when approaching the steady state, since we are considering $r = 10^{-3}$. {In particular, the pattern emerged from \textsc{Test-DD3} is in line with the results presented in \cite{astuto2023asymmetry}. On the other hand, the energy shows some oscillations, probably caused by the considered time step}.

In Fig.~\ref{fig_testAL}, we show the results with different values of the parameters $\alpha$ and $D^2$. From left to right we have the results of \textsc{Test-Al1} {with $\alpha = 0.75$ and $D^2 = 10^{-3}$,} \textsc{Test-Al2} {with $\alpha = 1.5$ and $D^2 = 5\cdot 10^{-4}$ and} \textsc{Test-DD3} {with $\alpha = 0.75$ and $D^2 = 5\cdot 10^{-4}$. We see more ramifications for $\alpha = 1.5$, but the respective energy plots show a more irregular behavior in time {since the oscillations are higher.}
	
	\begin{figure}
		\vspace{5mm}
		\centering
		\begin{minipage}{.3\textwidth}
			\centering
			\begin{overpic}[abs,width=0.9\textwidth,unit=1mm,scale=.25]{Figures_FEM/DD1.pdf}
				\put(0.5,0.1){\textcolor{white}{$\rightarrow\,x$}}
				\put(0.01,2){\textcolor{white}{$\uparrow$}}
				\put(0.01,7){\textcolor{white}{$y$}}
				\put(7,43){\textsc{Test-DD1}: $\|\mathbb{C}\|$}
			\end{overpic}
		\end{minipage}
		\begin{minipage}{.3\textwidth}   
			\centering
			\begin{overpic}[abs,width=0.9\textwidth,unit=1mm,scale=.25]{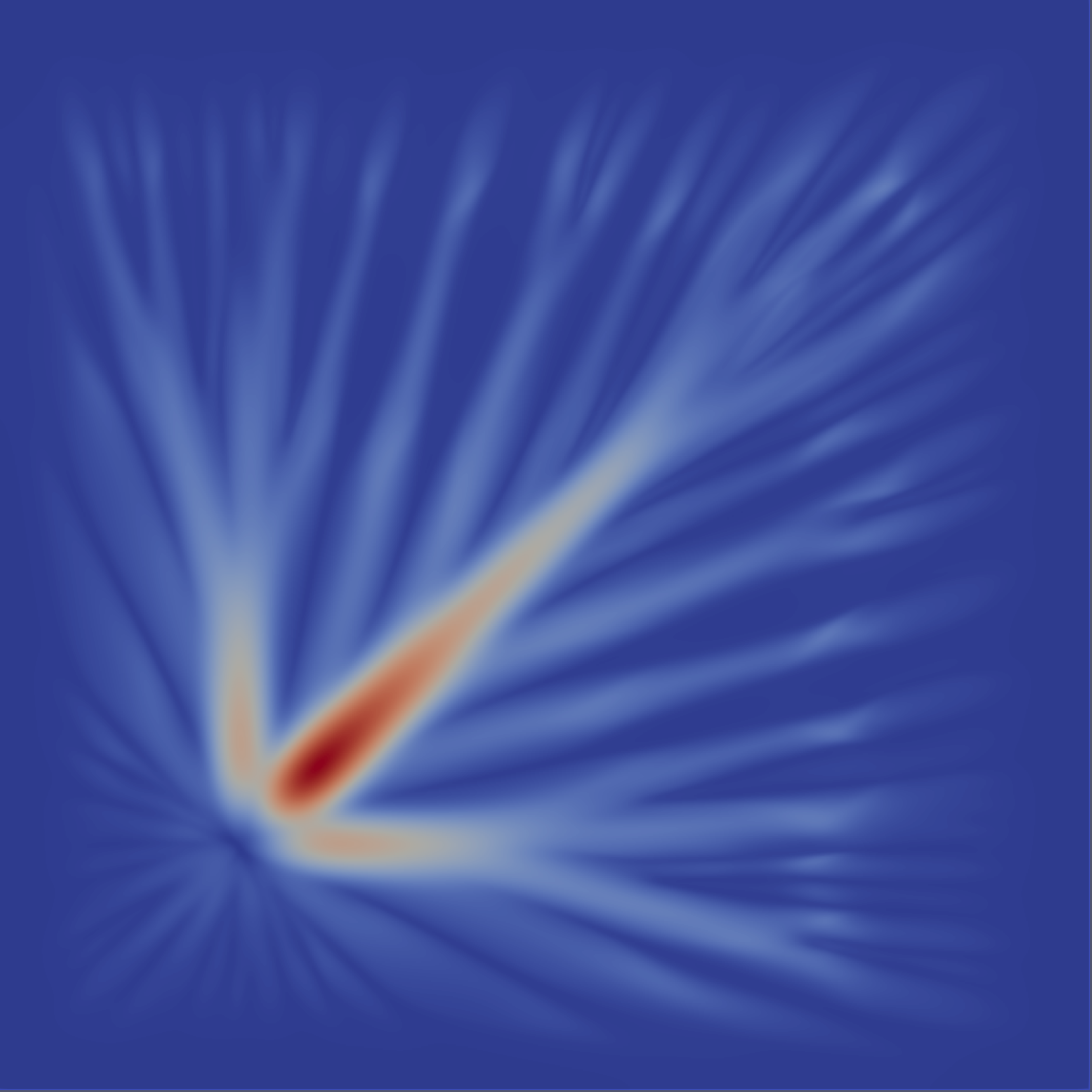}
				\put(0.5,0.1){\textcolor{white}{$\rightarrow\,x$}}
				\put(0.01,2){\textcolor{white}{$\uparrow$}}
				\put(0.01,7){\textcolor{white}{$y$}}
				\put(7,43){\textsc{Test-DD2}: $\|\mathbb{C}\|$}
			\end{overpic}
		\end{minipage}
		\begin{minipage}{.3\textwidth} 
			\centering
			\begin{overpic}[abs,width=0.9\textwidth,unit=1mm,scale=.25]{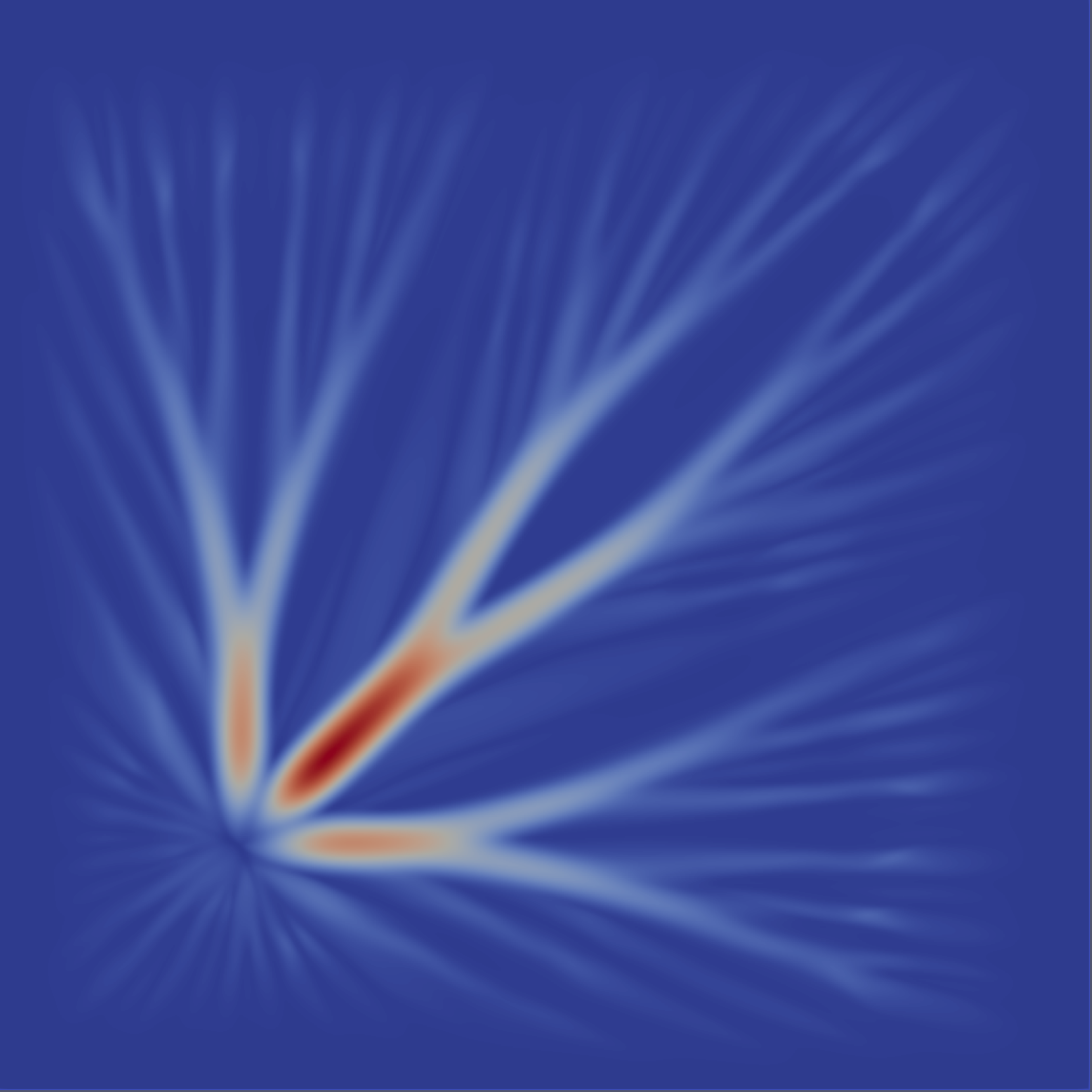}
				\put(0.5,0.1){\textcolor{white}{$\rightarrow\,x$}}
				\put(0.01,2){\textcolor{white}{$\uparrow$}}
				\put(0.01,7){\textcolor{white}{$y$}}
				\put(7,43){\textsc{Test-DD3}: $\|\mathbb{C}\|$}
			\end{overpic}
		\end{minipage}
		\begin{minipage}{.3\textwidth} \centering  \includegraphics[width=1.\textwidth]{Figures_FEM/DD1_energy.png}
		\end{minipage}
		\begin{minipage}{.3\textwidth}    \centering\includegraphics[width=1\textwidth]{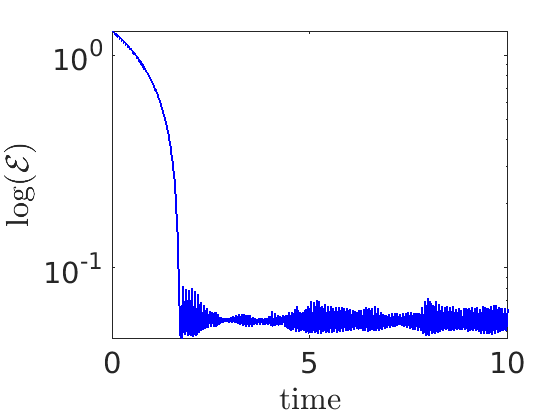}
		\end{minipage}
		\begin{minipage}{.3\textwidth}  \centering  \includegraphics[width=1\textwidth]{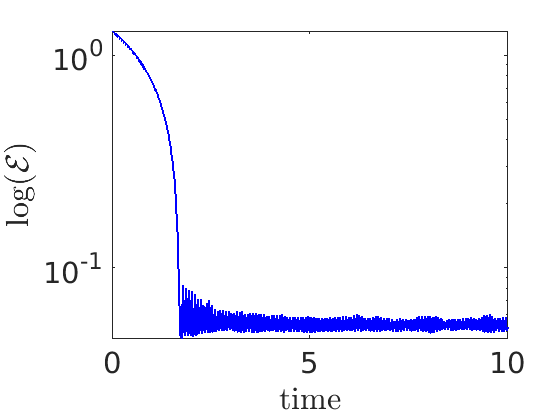}
		\end{minipage}
		\caption{\textit{In this figure we show three results considering different values of the squared diffusivity $D^2$. From left to right we have the results of }\textsc{Test-DD1} \textit{with $D^2 = 2\cdot 10^{-3}$,} \textsc{Test-DD2} \textit{with $D^2 = 10^{-3}$ and} \textsc{Test-DD3} \textit{with $D^2 = 5\cdot 10^{-4}$. Same format as Fig.~\ref{fig_testr}.}}
		\label{fig_testDD}
	\end{figure}
	
	\begin{figure}
		\vspace{5mm}
		\centering
		\begin{minipage}{.3\textwidth}   
			\centering
			\begin{overpic}[abs,width=0.9\textwidth,unit=1mm,scale=.25]{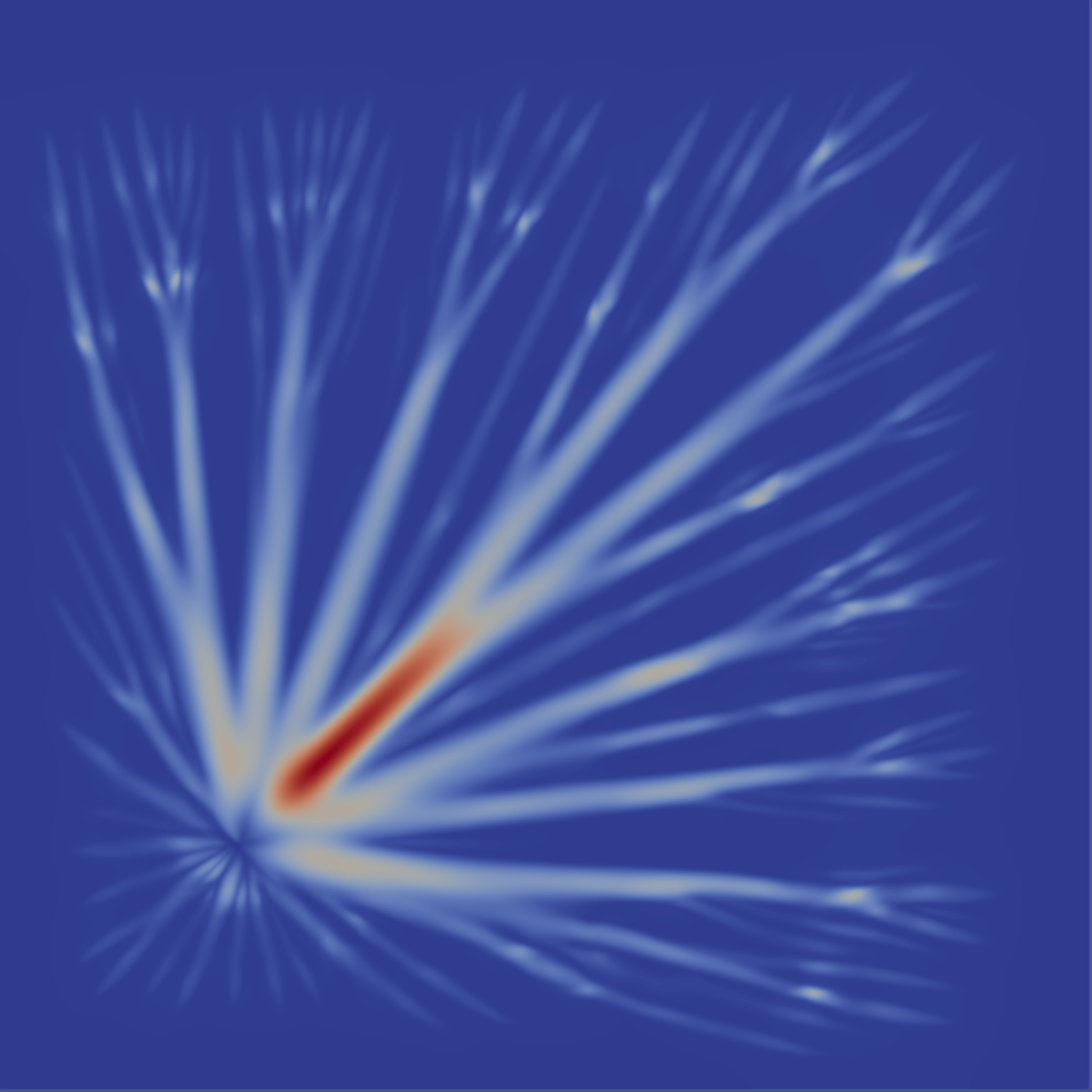}
				\put(0.5,0.1){\textcolor{white}{$\rightarrow\,x$}}
				\put(0.01,2){\textcolor{white}{$\uparrow$}}
				\put(0.01,7){\textcolor{white}{$y$}}
				\put(7,43){\textsc{Test-AL1}: $\|\mathbb{C}\|$}
			\end{overpic}
		\end{minipage}
		\begin{minipage}{.3\textwidth} 
			\centering
			\begin{overpic}[abs,width=0.9\textwidth,unit=1mm,scale=.25]{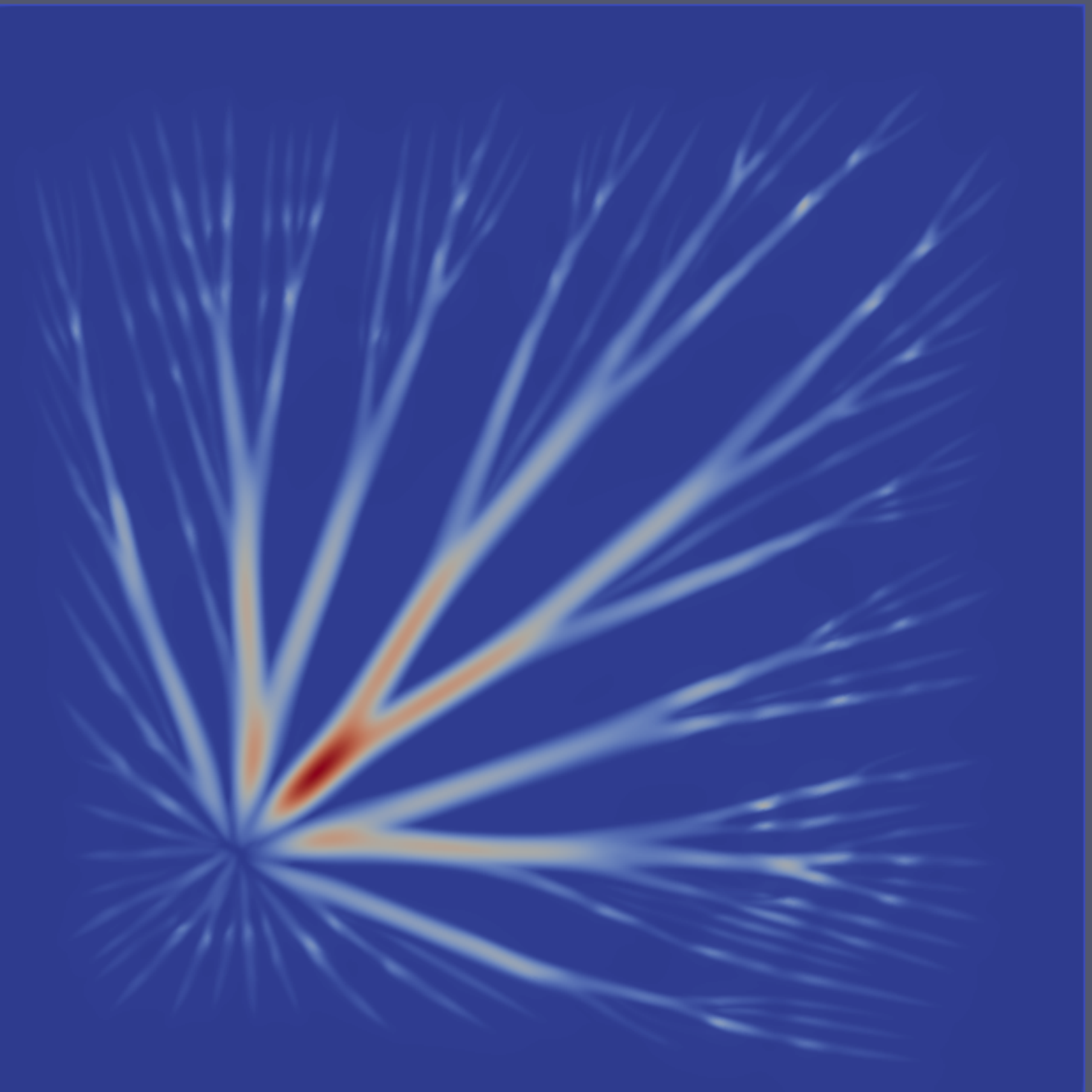}
				\put(0.5,0.1){\textcolor{white}{$\rightarrow\,x$}}
				\put(0.01,2){\textcolor{white}{$\uparrow$}}
				\put(0.01,7){\textcolor{white}{$y$}}
				\put(7,43){\textsc{Test-AL2}: $\|\mathbb{C}\|$}
			\end{overpic}
		\end{minipage}
		\begin{minipage}{.3\textwidth}
			\centering
			\begin{overpic}[abs,width=0.9\textwidth,unit=1mm,scale=.25]{Figures_FEM/DD3.pdf}
				\put(0.5,0.1){\textcolor{white}{$\rightarrow\,x$}}
				\put(0.01,2){\textcolor{white}{$\uparrow$}}
				\put(0.01,7){\textcolor{white}{$y$}}
				\put(7,43){\textsc{Test-DD3}: $\|\mathbb{C}\|$}
			\end{overpic}
		\end{minipage}    
		\begin{minipage}{.3\textwidth} \centering   \includegraphics[width=1\textwidth]{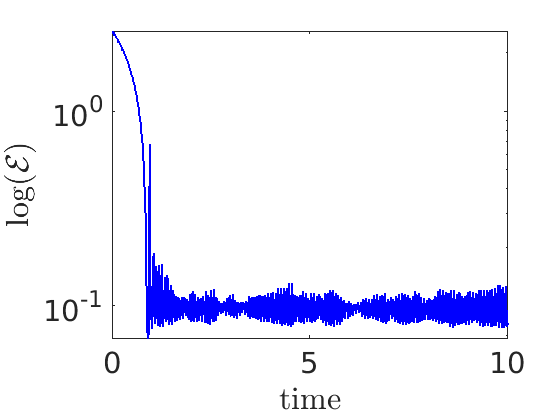}
		\end{minipage}
		\begin{minipage}{.3\textwidth} \centering   \includegraphics[width=1\textwidth]{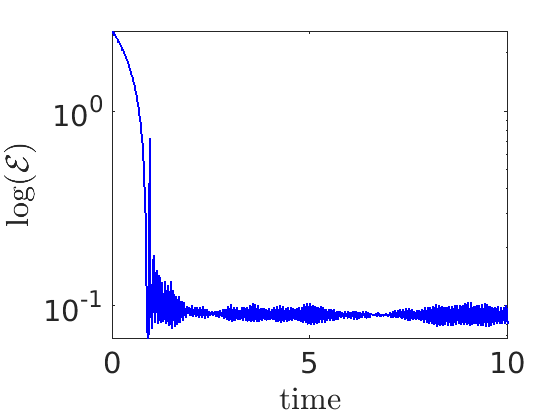}
		\end{minipage}
		\begin{minipage}{.3\textwidth} \centering  \includegraphics[width=1.\textwidth]{Figures_FEM/DD3_energy.png}
		\end{minipage}    
		\caption{\textit{In this figure we show three results considering different values for $\alpha$ and $D^2$. From left to right we have the results of} \textsc{Test-Al1} \textit{with $\alpha = 0.75$ and $D^2 = 10^{-3}$,} \textsc{Test-Al2} \textit{with $\alpha = 1.5$ and $D^2 = 5\cdot 10^{-4}$ and} \textsc{Test-DD3} \textit{with $\alpha = 0.75$ and $D^2 = 5\cdot 10^{-4}$. Same format as Fig.~\ref{fig_testr}.}}
		\label{fig_testAL}
	\end{figure}

	\section{Conclusions}
	\label{section_conclusion}
	In this paper we studied the Cai-Hu model, an elliptic-parabolic system to describe the formation of biological networks, and, in particular, of leaf venation networks. 
	
	We used a finite element method in space associated with a semi-implicit discretization in time, due to the nonlinearities of the model. The use of the library FreeFEM makes the presented space discretization very easy to implement. 
	The convergence rate in time and space is computed numerically and confirms the first order accuracy. 
	
	We studied the evolution of the energy, that decays in time, but it is non-monotonic when the parameter $r \to 0$. This depends on the non-solvability of the Poisson equation when $r=0$.  
	
	We remark that a big variety of parameters is involved in the system and each of them seems to influence the solution in different ways. A deeper investigation should focus on their roles in the equations.  
	
	Our future research will investigate how to improve the numerical method by computing the solutions in a monolithic way: for instance, this can be done by considering a mixed formulation. In this case, in order to treat the nonlinear terms in the system, an iterative method can be adopted, such as a Newton or a fixed point iterator.
	
	The ADI scheme, already applied in \cite{astuto2023asymmetry,astuto2022comparison} in its symmetric version, would improve the time accuracy of the scheme. {High order in time can be assured also with the use of implicit-explicit (IMEX) Runge-Kutta schemes. These methods are effective also in presence of stiffness in the system;} {in particular, the family of Stiffly Accurate (SA) schemes have been created for this purpose.}

Moreover, the use of mesh adaptivity could accurately track the conductivity flow by refining the mesh in those regions of interest where ramifications start appearing, {deeply examining the multiple scales hidden in the problem. These features are already under investigation, together with} {the design of} {a parallel solver able to accelerate} {the solution of the linear system at each time step.} {In this paper we were able to} {solve} {the equations for the three components of the tensor, plus one for the pressure, choosing a} {mesh size} {equal to 1/600.} {This resolution is already computationally expensive for a serial algorithm.}

\section*{Acknowledgments}
F. C. is partially supported by IMATI/CNR.


\end{document}